\documentclass[12pt]{amsart}

\makeatletter
\let\over\@@over
\makeatother
\let\cal\mathcal
\def\length{\operatorname{length}}

\usepackage{psfig}
\begin{document}

\title{On Two-Generator Satellite Knots}

\date{Version 1.0 of February 27, 1997} 

\author{Steven A. Bleiler}

\address{\hskip-\parindent Steven A. Bleiler, Mathematics Department,
Portland State University, Portland, OR  97207} 
\email{steve@math.pdx.edu}

\author{Amelia C. Jones}

\address{\hskip-\parindent Amelia C. Jones,  Mathematics Department,
Vassar College, Poughkeepsie, NY, 12601}
\email{amjones@vassar.edu}


\begin{abstract}
Techniques are introduced which determine the geometric
structure of non-simple two-generator $3$-manifolds from purely
algebraic data.  As an application, the satellite knots in the
$3$-sphere with a two-generator presentation in which at least one
generator is represented by a meridian for the knot are classified.   
\end{abstract}

\thanks{Research at MSRI is supported in part by NSF grant DMS-9022140}

\maketitle

\section{Introduction} 

We consider classical knot and
link exteriors as a context in which to introduce  techniques which
determine the geometric structure of non-simple $2$-generator
$3$-manifolds from purely algebraic data.  Here the application of
these  techniques are to a pair of related conjectures on the Heegaard
realization of generators.  These are the so called ``two generator is
tunnel one"  and ``meridional generator" conjectures. 

These two conjectures arise from two of the problems in R. Kirby's
famous 1978 list [K].   The first is a problem attributed to L. Moser
which asks for a geometric characterization of  those knots and links
whose fundamental group may be presented using two group elements. The
``two generator is tunnel one" conjecture purports to be  the solution
to this particular problem.  This conjecture is usually attributed to
M.  Scharlemann, who attributes it to A. Casson, who notes it is
implicit in a question asked in 1967 by F. Waldhausen [W]. The second is
a problem attributed to S. Cappell and J.  Shaneson which asks if
generators which can be represented (under a suitable choice of 
basepoint) by a simple closed curve on the boundary (the so-called {\it
meridional}  generators) correspond to the bridges of a bridge
decomposition of the knot exterior.  The  meridional generator
conjecture asserts an affirmative answer to Cappell and Shanneson's 
question although considered in the broader context of a generalized
notion of bridge  decomposition.  

Little is known on either of these conjectures. The only published
result on the first is due  to S. Bleiler [B] where the conjecture is
established for cable knots.  A single published  result exists for the
second as well.  This is due to M. Boileau and B. Zimmermann [BZ], who 
use Thurston's orbifold technology to show that a knot or link exterior
whose group is  generated by a pair of meridional elements is in fact
two bridge.

The principal result of this paper, Theorem 3.1, classifies those
satellites with a two  generator presentation in which one of the
generators is meridional.  Note that this includes  all known two
generator satellites.   The two generator is tunnel one and meridional 
generator conjecture for these knots follows as a corollary to the
classification.  It also  follows from this classification that a two
generator satellite knot out in this class is a  counter-example to the
two generator is tunnel one conjecture.  The key to the classification 
is an algebraic theorem given in 3.10.  

The paper is organized as follows.  The first and second sections
respectively contain the  topological and algebraic backgound necessary
for the classification theorem which is stated and proved in section
three. 

The techniques developed here apply more generally to the
classification of non-simple
$3$-manifolds.  The details of this will appear in a subsequent
publication.

\section{Topological Background} 

\medbreak\noindent
{\bf 1.1.}  Begin by recalling B. Clark's construction of a knot whose 
fundamental group is at most $g$-generator [C].  Start with a genus $g$
Heegaard  splitting of the $3$-sphere and consider one of the
handlebodies as the union of a single 
$0$-handle and $g$ $1$-handles.  Form the exterior of a knot by removing
the $0$-handle and one of the $1$-handles.  The cores of the remaining
$1$-handles form a set  of {\it unknotting tunnels } for the knot $K$. 
The minimal cardinality of a set of  unknotting tunnels for $K$ is
called the {\it tunnel number} of $K$.  It follows that the  rank of the
fundamental group of a  tunnel number $n$ knot is at most $n+1$.  It is
a consequence of  Dehn's lemma that the rank of the group of a tunnel
number one knot is exactly two. Clark's construction can be generalized
to form complements of $n$ component links by decomposing the genus
$g$  handlebody into a collection of $n$ $0$-handles and $g+n-1$
$1$-handles and deleting a collection of $n$ $0$-handles and $n$
$1$-handles. It is  occasionally useful to consider this construction
in the context of general 3-manifolds.     

\medbreak\noindent
{\bf 1.2.}  Next, recall H. Doll's generalization of Schubert's notion
of bridge  number [D].  A properly embedded arc $\alpha$ in a
handlebody is said to be {\it  trivial } if there exists an arc $\beta$
in the boundary of the handlebody whose endpoints  agree with those of
$\alpha$ such that $\alpha \cup \beta$ bounds a disc in the
handlebody.   A knot in a $3$-manifold $M$ has a {\it
$(g,b)$-decomposition} or alternately, is said to be in {\it
$b$-bridge  position with respect to a Heegaard surface $F$  of genus
$g$} if the knot intersects the  closure of each complementary
handlebody to $F$  in $b$ trivial arcs.  One now defines  the {\it genus
$g$ bridge number} for a link $L$ to be the minimal number $b$ for
which  the link has a $(g,b)$ decomposition with respect to a genus $g$
splitting of the ambient  manifold $M$.  Note that Schubert's original
definition of bridge number is the above definition applied to $S^3$
with $g = 0$.  Also note that when $b \geq 1$, by stabilizing the
Heegaard splitting one can  always change bridges to genus, that is, a
link with a $(g,b)$ decomposition always has a $(g+1, b-1)$
decomposition.  That the converse is false is demonstrated by the $(-2,
3, 7)$  pretzel knot which has no $(0,2)$ decomposition (i.e. is not two
bridge), but as the closure  of the 7-string braid
$(\sigma_1\sigma_2)^{-1}\cdot (\sigma_1\sigma_2 \sigma_3 \sigma_4
\sigma_5 \sigma_6)^3$ has a $(1,1)$ decomposition.  See Figure 1.  Also
note  from Figure 1 that as a pretzel knot, the $(-2,3,7)$ pretzel also
has a
$(0,3)$  decomposition.

\begin{figure}
\centerline{%
\psfig{figure=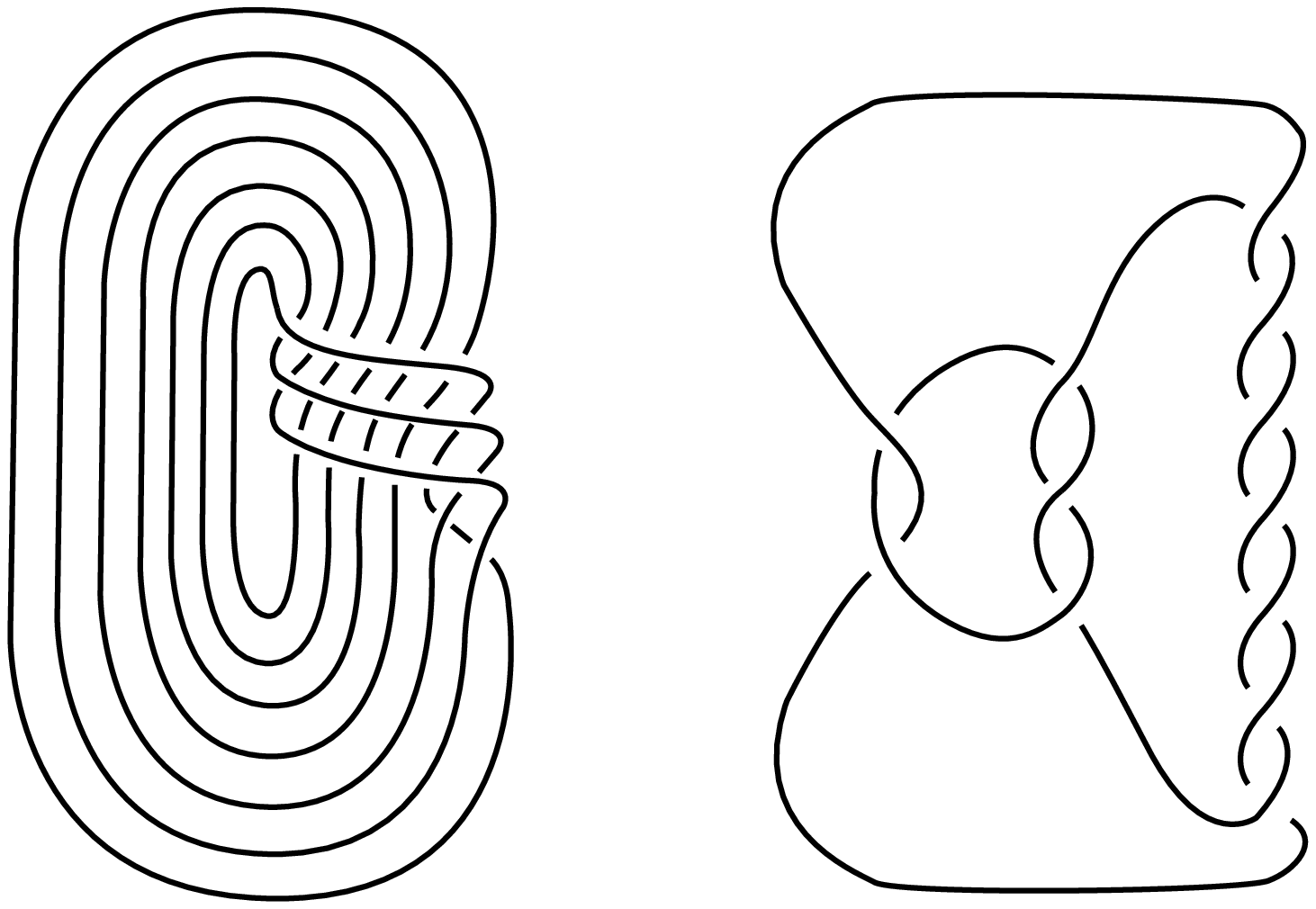,width=4in}%
}  
\caption{}
\end{figure}

\noindent For convenience, we refer to a knot with a $(1,1)$
decomposition as {\it one bridge with  respect to a torus}.  Note that
such knots are naturally tunnel one, see [G].  

\medbreak\noindent
{\bf 1.3.}  We now develop the terminology that describes how these
geometric  concepts appear in the  algebraic data.  An element of the
fundamental group of a link  exterior is said to be {\it meridional} if
there is a choice of basepoint so that the element  is represented by a
simple closed curve in the boundary.  Note that an element can be 
meridional and yet not lie in the fundamental group of the boundary of
the knot, as  exemplified by the two bridge knot and link exteriors. The
fundamental group of a two  bridge link can be generated by two
meridional elements, both of which cannot lie in the fundamental group
of the boundary.       

A link exterior is said to have a {\it $(g,b)$ presentation} if its
fundamental group has a  presentation with $g+b$ generators, $b$ of
which are meridional.  Note that if $b \geq 1$  that a knot with a
$(g,b)$ decomposition has a $(g,b)$ presentation.   The converse is:

\medskip
\noindent{\bf Meridional Generator Conjecture.} {\it A knot with a
$(g,b)$ presentation has a 
$(g,b)$ decomposition.}
\medskip

Remark:  That a knot with a $(g,0)$ decomposition need not have a
$(g,0)$ presentation is  exemplified by the torus knots.   However, as
of this writing, no counterexample is known  to the conjecture that a
knot with a $(g,0)$ presentation has a $(g,0)$ decomposition. 

Remark:  Cappell and Shaneson's original problem that n-meridional
generators implies 
$n$-bridge is just the special case $g = 0$.

\medbreak\noindent
{\bf 1.4.}  The tunnel one satellite knots in $S^3$  were classified in
1988 by  K.Morimoto and H. Sakuma in [MS] and the tunnel one satellite
knots in general $3$-manifolds were classified in 1995 by J. Neil [N]. 
Per Morimoto and Sakuma's  classification, a satellite knot $K$ is
tunnel one when it satisfies the following three conditions.  First,
the exterior of $K$ must be the union of a simple link  exterior $X_A$
and a torus knot exterior $X_B$.  Second, the link exterior $X_A$ must
be  attached to $X_B$ in a manner that takes a meridian of $X_A$ to the
Seifert fibre of $X_B$.  Finally, the link  exterior $X_A$ must embed
in $S^3$ as the exterior of a $2$-bridge link other than the un- or 
Hopf links.  As can be seen from the construction, all the knots
obtained in  this manner are one bridge with respect to the torus knot's
torus.  See Figure 2 where  this is illustrated for the Rolfsen-Bailey
knot, i.e. the one obtained by attaching the  exterior of the two bridge
link
$4 \over 1$ to the exterior of the left hand trefoil, as  described
above.  It is interesting to note that M. Eudave-Munoz was able to
recover this  classification by first proving that a tunnel one
satellite in $S^3$ has a $(1,1)$-decomposition, see [E].

\begin{figure}
\centerline{%
\psfig{figure=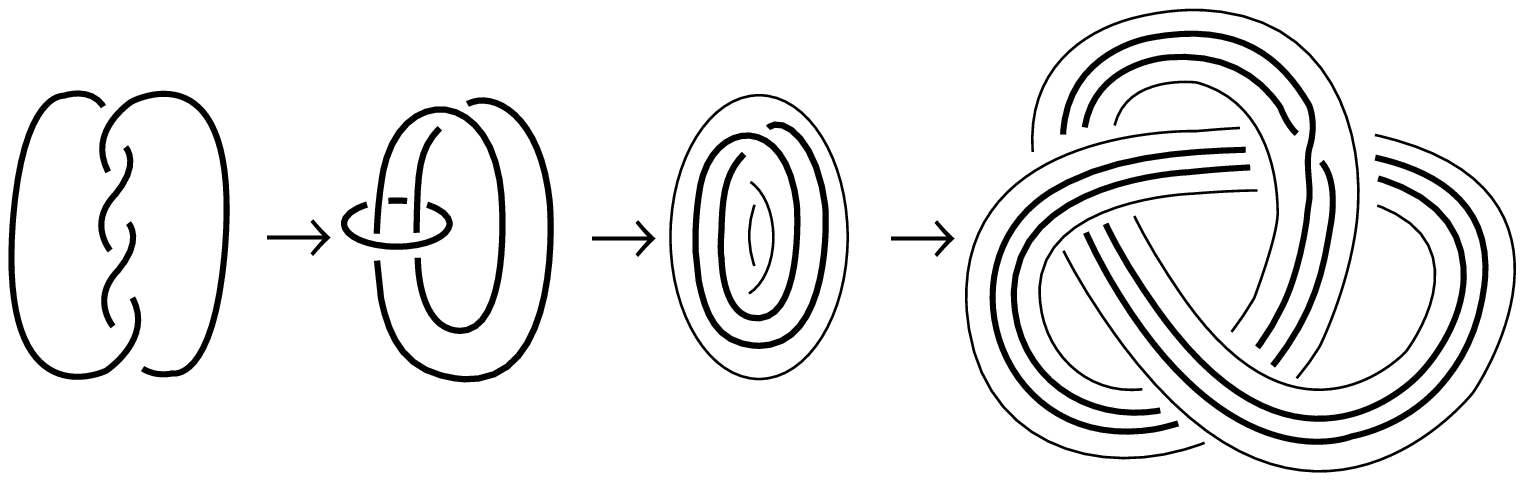,width=5in}%
}  
\caption{}
\end{figure}

\section{2. Algebraic Background} 

\medbreak\noindent
{\bf 2.1.}    Here we review some of the properties  of groups needed in
future discussion.  Let $A$ denote a group.  A subgroup  $S$ of $A$ is
called {\it proper} if $S$ is nontrivial and there exist elements  of
$A$ outside of $S$.  A proper subgroup $S$ is said to be {\it malnormal}
in $A$ if the intersection of subgroups $aSa^{-1}$ and $S$ is trivial if
and only  if $a$ is an element of $A$ outside of $S$.  Thus $S$
intersects each of its nontrivial conjugates trivially.  An element $a
\in A$ is said to be {\it central} if $a$ commutes with each element of
$A$.  A {\it torsion} element of $A$ is an element $a$ satisfying the
relation $a^n = 1$ for some nonzero integer $n$.    

\medbreak\noindent
{\bf 2.2.}  The {\it {rank}} of a finitely presentable group is the
minimal number of generators needed to present the group.  A group $A$
is called an {\it {n-generator}} {\it{group}} if $A$ is known to have
rank $n$.  The groups of rank one are the finite cyclic groups together
with the infinite cyclic group.  The groups of higher rank form a much
more complicated class.   Indeed, determining the rank of a group is in
general a difficult problem.  

\medbreak\noindent
{\bf 2.3.}   The following concepts are presented in greater detail in
Magnus,  Karrass, and Solitar [MKS] and are recalled here for the
convenience of the reader.  Let $A*_CB$ denote the free product with
amalgamation of the non-trivial groups $A$ and $B$ along the group $C$,
that is, $C$ is isomorphic to a proper subgroup $A_C$ of  $A$ and also
isomorphic to a proper subgroup $B_C$ of $B$.  The subgroups $A_C$ and
$B_C$ are identified with the group $C$ via the respective
isomorphisms.   When the context is clear, we simply refer to the
subgroup
$A_C$ of $A$ by $C$, and similarly the subgroup $B_C$ of $B$ will be
referred to by $C$.  Note that group $G$ contains each of $A$, $B$ and
$C$ as proper subgroups. The groups $A$ and
$B$  are called the {\it free} {\it factors} for the group $A*_CB$ and
the group $C$ is called the {\it amalgamated} {\it subgroup}.  

Any element $x$ of $A*_CB$ may be uniquely expressed in {\it normal}
{\it form} as 
$x=x_1x_2\cdots x_mc$ where $c$ is an element of $C$ and the $x_i$ are 
alternately elements of a fixed set $T_A$ of right transversals for
the  (right) cosets of $C$ in $A$ and a corresponding fixed set $T_B$ of
right  transversals for the (right) cosets of $C$ in
$B$.  The coset $C$ in  either of $A$ or $B$ always receives the
transversal represented by 1 and the transversal representative 1 never
appears in the normal form word for  $x$.  If $x$ is expressed in the
above normal form, then $x$ is said to have  length $m$, and one writes 
$\length(x)=m$.  If $x_1$ is an element of $A$,  then $x$ is said to
begin in $A$.  If $x_m$ is an element of $A$, then $x$  is said to end
in $A$.   Similarly, if $x_i$ is an element of
$B$ for $i=1$ or $i=m$, then $x$ is said to begin or end respectively in
$B$.  If $x$ begins in $A$ and ends in $B$, then $x$ can be written 
concretely in the normal form $x=a_1b_2\cdots b_mc$.  If $x$ can be
written in this form, then $x$ is said to be a word of length $m$
beginning in $A$ and ending in  $B$, and similarly for the other cases.  

Two words $x$ and $y$ written in normal form are multiplied by
concatenation.  The product is subsequently converted to normal form
by  moving the elements of $C$ to the right through the word $y$.  If
the word  $x=x_1\cdots x_nc$ ends in the same group that the word
$y=y_1\cdots  y_mc'$ begins in, then the product $x_ncy_1$ is an
element of one of $A$ or 
$B$.  If $x_ncy_1$ is an element of the complement of $C$ in that
factor,  then $xy$ is said to have  {\it {amalgamation}} and
$\length(xy)= n+m-1$.  If  the product $x_ncy_1 = \bar c$, an element of
$C$, then the word $xy$ has  {\it {cancellation}} and the product
$x_{n-1} \bar c y_2$ is then examined  for additional cancellation
and/or amalgamation and  $\length(xy)\leq m+n- 2$.  If $x$ ends in a
different group than $y$ begins in, there is no  amalgamation or
cancellation in $xy$ and $\length(xy)=m+n$. 
  
Now recall the results of [BJ] which will be relevent to the discussion
at hand.

\medbreak\noindent
{\bf 2.4.  Lemma [BJ, 1.5].} {\it  Suppose that the group $A*_CB$ is
generated by
$n$
 elements $g_1$, $g_2$, \dots $g_n$.   The set $S_A$ consisting of
those  transversals of $T_A$ appearing in any of the $g_i$ together
with all of  the elements of $C$ generates the group $A$. }

\medbreak\noindent
{\bf 2.5.}   By symmetry, the group $B$ may be generated by the set
$S_B$ consisting of all transversal elements of $T_B$  that appear in
at least one of the generators $g_i$ together with the elements of $C$. 

\medbreak\noindent
{\bf 2.6.  Theorem [BJ 3.1].}  {\it Let $G = A*_C B$  be the free
product with amalgamation of the torsion free group $A$ and the group
$B$.  If the amalgamated subgroup $C$ is of rank at least two and
malnormal in $A$ and if $G$ has rank two, then either $B$ contains
torsion elements or
$G$ has a generating pair of one of four types (listed here as in {\rm
[BJ, 5.1]})}:

$p_1:  \lbrace bc, ac_0 \rbrace$

$p_2:   \lbrace c, a_1b_2c_0 \rbrace$

$p_3:  \lbrace bc, a_1b_2c_0 \rbrace$

$p_4:  \lbrace bc, a_1b_2a_3c_0 \rbrace$

\medbreak\noindent
{\bf 2.7.  Proposition [BJ 5.1.1].}  {\it If $G$ has a generating pair
of type $p_1$ and
$C$ is abelian, then there is a minimal positive integer $m$ such that
the element $(bc)^m$ is an element $\tilde c$ of $C$ which is central in
$B$.}

\section{The $(1,1)$ Presented Satellite Knots} 

\medbreak\noindent
{\bf 3.1.  Theorem.} {\it The exterior of a $(1,1)$ presented satellite
knot in the
$3$-sphere decomposes as the union along a torus of a simple link
exterior $X_A$ and a  torus knot exterior $X_B$.  Moreover, the link
exterior $X_A$ embeds in the $3$-sphere as the complement of a two
bridge link and is attached to $X_B$ so that a meridian for this two
bridge link is identified with the Seifert fibre of the torus knot
exterior.}

Combining 1.4 and 3.1 gives the following corollories.

\medbreak\noindent
{\bf 3.2.  Corollary.}  {\it A $(1,1)$ presented satellite knot in the
$3$-sphere is one bridge with respect to a torus, in particular, is
tunnel one.}

\medbreak\noindent
{\bf 3.3.  Corollary.}  {\it A two-generator satellite knot in the
$3$-sphere with no
$(1,1)$ presentation has tunnel number at least two.}

\medbreak\noindent
{\bf 3.4.}  Let $K$ be a $(1,1)$ presented satellite knot in $S^3$.  If
$K$ is a cable knot, then by [B], $K$ is the $(spq \pm 1,s)$ cable on
the $(p,q)$ torus knot. Here, the space $X_A$ embeds in $S^3$ as the
complement of the two bridge link corresponding to the  even integer $2s
\over 1$.  The conclusion of 3.1 follows.  

Without loss of generality, the knot  $K$ is not cabled.  Denote the
exterior of $K$ by $X_K$.  The space $X_K$ has a decomposition by
essential tori into simple pieces, see [JS] or [J].  One of these simple
pieces has $\partial X_K$ as one of its boundary components.  Denote
this piece by $X_A'$.  Thurston's geometrization for Haken manifolds
[T] implies that $X_A'$ is either a Seifert fibred space or its
interior admits a hyperbolic metric.  If not hyperbolic, the space
$X_A'$ contains a non-boundary parallel essential annulus, vertical
with respect to the Seifert fibration on
$X_A'$, with both boundary components lying on the boundary of $X_K$. 
As $K$ is not a cabled knot, it follows that the knot $K$ is composite,
contradicting Norwood's result [N] that two-generator knots are prime. 
Hence the space $X_A'$ is hyperbolic.  

The space $X_A'$ has at least two boundary tori.  Let the torus $T$ be a
boundary torus that is different from $\partial X_K$.  This torus $T$
separates the space $X_K$ into two pieces.  Let $X_A$ be the component
of $X_K - T$ containing $X_A'$ and $X_B$ the complementary component. 
Since T is a separating torus, by the Seifert-Van Kampen theorem, the
group $G=\pi_1(X_K)$ splits as the free product $A *_C B$ where $A$ is
the group $\pi_1(X_A)$,
$B$ is the group  $\pi_1(X_B)$ and $C$ is the group of the torus
$T$.       

\medbreak\noindent
{\bf 3.5.  Proposition.} {\it The group $C$ is a malnormal subgroup of
$A$.}

\medbreak\noindent
{\bf 3.6.}  ({\it Proof of 3.5.})  Let $g$ be an element of the group
$A$ such that
$gcg^{-1}=c'$ for  elements $c,c' \in C$.  Hence there exists a proper,
$\pi_1$-injective map
$f$ of the  annulus $S^1 \times I$ into $X_A$.  Put the image of $f$
into general position with respect to the  decomposing tori for $X_A$
inherited from the torus decomposition of
$X_K$.  While the  intersection of the image of $f$ with one of these
tori may be hideous on the torus, the  intersection on the annulus may,
by standard techniques, be reduced to a family $\cal F$ of  parallel
essential simple closed curves on the annulus.  The restricition of $f$
to a  complementary component of $\cal F$ which contains a component of
$\partial (S^1
\times  I)$ is a proper, $\pi_1$-injective map of an annulus into the
hyperbolic manifold
$X_A'$,  hence the image of this restricted $f$ can be isotoped into the
torus $T$.  Induction then shows that the entire image $f$ can be
isotoped into $T$ and hence that $g \in C$.
\hfill $\diamondsuit$
 
\medbreak\noindent
{\bf 3.7. Proposition.}  {\it The group $G= \pi_1(X_K)$ has a
generating pair of type $\lbrace ac_0, bc \rbrace$.}

\medbreak\noindent
{\bf 3.8.}  ({\it Proof of 3.7.})  As the group $G$ is
$(1,1)$-presented, the group $G=A*_C  B$ has a generating pair $\lbrace
g_1, g_2 \rbrace$.  Without loss of generality, the generator $g_1$ is
meridional and hence it has normal form $ac$.  Thus there are exactly
five possibilites for  the normal forms of the generating set $\lbrace
g_1, g_2 \rbrace$.  These are:

$p_1:  \lbrace ac, bc \rbrace$

$q_2:   \lbrace ac_0, a_1b_2 \cdots b_{2k}a_{2k+1}c \rbrace$

$q_3:  \lbrace ac_0, a_1b_2 \cdots b_{2k}c  \rbrace$

$q_4:  \lbrace ac_0, b_1a_2 \cdots a_{2k}c  \rbrace$

$q_5:  \lbrace ac_0, b_1a_2 \cdots a_{2k}b_{2k+1}c  \rbrace$

By inverting the second element, we can consider a pair of type $q_4$ to
be of type $q_3$.   But by [BJ 4.4], a pair of type $q_3$ cannot
generate $G$ and by [BJ 4.5] neither can a  pair of type $q_5$.  So if
$\lbrace g_1, g_2 \rbrace$ is not of type $p_1$, it must be of type
$q_2$.  Note that by an argument similar to that in 3.6, the element
$a_1^{-1}ac_0a_1$ is never an element of the  subgroup $C$.  This is
essentially the argument that a hyperbolic manifold contains no
essential annuli. Thus conjugating a pair of type $q_2$ by the element
$a_1^{-1}$ produces a generating set of type $q_4$  or $q_5$.  Thus
$\lbrace g_1, g_2
\rbrace$ is not of type $q_2$.  The desired result follows.
\hfill $\diamondsuit$

\medbreak\noindent
{\bf 3.9.}  From 3.7 and 2.7, the element $\tilde c = (bc)^m$ is central
in $B$.  It follows from G. Burde and H. Zieschang [BZ] that $X_B$ is
the complement of a torus knot and that the element $\tilde c$ is a
power of the class of the torus knot's regular Seifert fibre.  There is
a choice  of basis $\lbrace x, y \rbrace$ for the group $C$ so that $x$
is the class of a meridian for the torus knot and $y$ coincides with
$\tilde c$ if $\tilde c$ is a primitive element of $C$ and is a
primitive root of $\tilde c$ otherwise.  

\medbreak\noindent
{\bf 3.10. Theorem.}  {\it The elements $ac_0$ and $\tilde c$ generate
the group
$A$.  Hence $A$ has a $(0,2)$ presentation.  In particular, the space
$X_A$ embeds in $S^3$ as the exterior of a two bridge link with
meridians representing $ac_0$ and a primitive root of $\tilde c$.}

\medbreak\noindent
{\bf 3.11. Lemma.}  {\it The elements $ac_0$ and $\tilde c$ generate a
subgroup of
$A$ containing an element of the form $x^q$ for some nonzero integer
$q$.}

\medbreak\noindent
{\bf 3.12.}  ({\it  Proof of 3.11.})  Assume not.  For any integer $p$,
an element of form
$x^p$  is an element of $C$ outside of the cyclic subgroup generated by
$\tilde c$, and as an element of $G$ can  be expressed as a word $w$ in
the generating pair $\lbrace ac_0, bc \rbrace$.  The word $w$ can be
written in the following form:

$w = (bc)^{p_1}(ac_0)^{p_2}(bc)^{p_3} \cdots
(ac_0)^{p_{2j}}(bc)^{p_{2j+1}}$.  

\noindent The integers $p_1$ and $p_{2j+1}$ are allowed to be zero, but
all other integers
$p_i$ are nonzero.  As $x^p$ is an element of $C$ the word $w$ admits
amalgamations and cancellations that reduce its length to zero.  We now
begin a series of rewritings of the word $w$.  Note that if $p_{2i+1}$
is a multiple of the integer $m$, then we may replace the phrase
$(bc)^{p_{2i+1}}$ by a power of $\tilde c$.  Do this throughout the
word $w$ and combine maximal length  phrases within $w$ that are words
in the elements $ac_0$ and $\tilde c$.  The word $w$ now has the
following appearance:

$w = (bc)^{p_1}w_2(bc)^{p_3} \cdots w_{2l}(bc)^{p_{2l+1}}$

\noindent Each word $w_i$ is a word in the elements $ac_0$ and $\tilde c
= (bc)^m$.  With the exception that $p_1$ and $p_{2l+1}$ are allowed to
be zero, no power $p_{2i+1}$ is a multiple of the integer $m$, else it
is combined into a neighboring $w_{2i}$ or $w_{2i+2}$. If $w$ has form
$(bc)^0w_2$, then we have established 3.11.  Hence we assume that each
phrase 
$w_{2i}$, if not in $A-C$ is a power of $\tilde c$.  Since $\tilde c$
and $bc$ commute, all pwers of $bc$ my be moved right through the word
$w$.  If
$w_{2i}$ is a power of $\tilde c$, any occurence of the phrase
$(bc)^{p_{2i-1}}w_{2i}(bc)^{p_{2i+1}}$ may be rewritten as
$w_{2i}(bc)^{(p_{2i-1}+p_{2i+1})}$, and then reindexed as necessary. 
In this way, the word
$w$ has been rewritten in the form:

$w = (bc)^{p_1}w_2(bc)^{p_3} \cdots w_{2n}(bc)^{p_{2n+1}}$

\noindent Each word $w_i$ is a word in $A-C$.  With the exception that
$p_1$ is allowed to be zero and $p_{2l+1}$ could be  an arbitrary
integer, no other power $p_{2i+1}$ is a multiple of the integer $m$. 
Notice that $w$ has now been written in a form admitting no further
cancellation or amalgamation.  Thus if $w$ is to have length zero, it
must have form $(bc)^p$ where $p$ is a multiple of $m$ or $(bc)^0w_2$. 
But either of these options produce a power of $\tilde c$, a
contradiction.
\hfill $\diamondsuit$

\medbreak\noindent
{\bf 3.13.}  ({\it Proof of 3.10.})  By 3.11, the elements $ac_0$ and
$\tilde c$ generate an element of form $x^p \in C$.  If $|p|=1$, then by
2.4, the elements  $ac_0$ and $\tilde c$  generate the entire group
$A$.  Assume that $ac_0$ and $\tilde c$ do not generate the element $x$
and let $k$ be the minimal positive integer for which the elements
$ac_0$ and
$\tilde c$  generate the element $x^k$.  It follows from the argument in
3.12 that since the pair $\lbrace ac_0,  bc, \rbrace$ generate the group
$G$, that the elements $bc$ and $x^k$ must generate an element of form
$x^j$ where $1 \leq j \leq k-1$.  Thus, Theorem 3.10 follows from the 
proposition below.

\medbreak\noindent
{\bf 3.14. Proposition.}  {\it The elements $bc$ and $x^k$ cannot
generate an element of form $x^j$ where $1 \leq j \leq k-1$.}

\medbreak\noindent
{\bf 3.15.}  ({\it Proof of 3.14.}) As noted in 3.9, $B$ is the group of
a torus knot.  Hence the group $B$ decomposes as the free product with
amalgamation of form $E *_D F$, where  each of the groups $D$, $E$ and
$F$ are infinite cyclic.   Geometrically, this corresponds to the
familiar picture of constructing a torus knot exterior by glueing
together two solid tori along an annulus in the boundary of each.  Each
infinite cyclic group $E$ and $F$ corresponds to the group of a solid
torus.  Let $h:B \longrightarrow {\bf Z}_p * {\bf Z}_q$ be the
homomorphism defined by quotienting the group $B$ by the infinite
cyclic subgroup generated by the class of the regular Seifert fibre. By
abusive  notation we both
$e_i$ and $h(ei)$ will be denoted by the symbol $e_i$. Similarly, both
$f_i$ and $h(fi)$ will be denoted by the symbol $f_i$, the context
making the meaning clear.  Note  for later use that as $p$ and $q$ are
the indices of the exceptional fibres in the Seifert fibration of a
torus knot exterior, they are relatively prime. Thus one of these
integers must be odd and so one of ${\bf Z}_p$ or ${\bf Z}_q$ has no
element of order two.

As the meridian for the torus knot, the element $x$ has normal form
$e_1f_2d \in E*_D F$, and so the  element $h(x)$ has normal form
$e_1f_2$ in the free product ${\bf Z}_p * {\bf Z}_q$. Note that  the
normal form for the element $x^k$ is then $(e_1f_2)^k$, in particular
has length $2k$.  Since $(bc)^m = \tilde c$, the element $h(bc)$ lies
in either a free factor of
${\bf Z}_p * {\bf  Z}_q$ or a conjugate thereof.  Each of these cases
is examined in turn.  

\medbreak\noindent
{\bf 3.15.1.}  Suppose the element $h(bc)$ lies in a free factor of
${\bf Z}_p * {\bf Z}_q$.  Without loss of generality, $h(bc) = e'$. 
Let $w$ be a word in the elements $h(bc)$ and
$h(x^k)$.  Then  we may assume that $w$ appears as follows:

$w = e'{}^{p_1}h(x^k)^{p_2}e'{}^{p_3} \cdots
h(x^k)^{p_{2j}}e'{}^{p_{2j+1}}$. 

\noindent With the possible exception of $p_1$ and $p_{2j+1}$, all
integers $p_i$ are nonzero and such that no power of $e'$ is trivial. 
If $w$ is a nontrivial power of the element
$e'$, note that $\length(w)$ in ${\bf Z}_p * {\bf Z}_q$ is one.  If $w$ is
a power of $h(x^k)$, say $h(x^k)^p$, then $\length(w) = |2kp|$.  We now
assume that $w$ has neither of these forms and compute the length of
$w$.  Before computing lengths, some preliminaries  are required:

Given a fixed word $w$ as above, associate a series of integers to $w$
as in [BJ].  In  particular, let  $P$ denote the associated sequence
$\lbrace p_2, p_4, \cdots p_{2k} 
\rbrace$.  The integer $\sigma $ denotes the number of times  that the
sequence $P$ changes algebraic sign.  For example, if $w$ is the  word
$e_1^8h(x^k)^4e_1^{-2}h(x^k)^{227}e_1h(x^k)^{- 88}e_1h(x^k)$, then the
sequence  $P$ is
$\lbrace 4, 227, -88, 1 \rbrace$ and the associated integer $\sigma$ is
$2$.  Note that the relation $0 \leq  \sigma \leq j-1$ holds.  Further,
for  each integer $2i + 1$ satisfying $0
\leq i \leq j$, let $\epsilon_{2i+1}$ denote  the length of the phrase
$e_1^{p_{2i+1}}$.  Also note that with the possible exception of $p_1$
and $p_{2j+1}$ all these integers are one.  Without amalgamation or
cancellation, the  ``uncancelled" length of the word $w$ is given by:

\centerline{$ \sum_{i=1}^{i=j} |p_{2i}|2k + \sum_{i=0}^{i=j} 
\epsilon_{2i+1} \geq 2kj + (j-1) + \epsilon_1 + \epsilon_{2j+1}$ }

In what follows, assume that the integer $k$ is at least two.  Note that
as the sequence $P$ proceeds from $p_{2i}$ to $p_{2(i+1)}$, the phrase
$h(x^k)^{p_{2i}}e'{}^p_{2i+1}h(x^k)^{p_{2i+2}}$  may admit cancellation
or amalgamation as outlined below:

\smallbreak\noindent
{\bf $\cdot$} If both $p_{2i}$ and $p_{2(i+1)}$ are positive,
the phrase
$h(x^k)e'{}^p_{2i+1}h(x^k) = (e_1f_2)^ke'{}^p_{2i+1}(e_1f_2)^k$ may have
cancellations and amalgamation.  By 3.15, there are at most two
cancellations and a single amalgamation yielding a maximal length
reduction of five.  

\smallbreak\noindent $\cdot$ If both $p_{2i}$ and $p_{2(i+1)}$ are negative, then
again by 3.15, the phrase $$h(x^k)^{-1}e'{}^p_{2i+1}h(x^k)^{-1} =
(f_1e_2)^{-k}e'{}^p_{2i+1} (f_1e_2)^{-k}$$ has at most two cancellations
and a single amalgamation yielding a maximal length reduction of five.

\smallbreak\noindent $\cdot$ If $p_{2i}$ is positive and $p_{2(i+1)}$ negative, then
writing out the phrase $h(x^k)e'{}^p_{2i+1}h(x^k)^{-1} =
(e_1f_2)^ke'{}^p_{2i+1}(e_1f_2)^{-k}$ shows there is no cancellation or
amalgamation.

\smallbreak\noindent $\cdot$ If $p_{2i}$ is negative and $p_{2(i+1)}$ positive,
then since the group
${\bf Z}_p$ is abelian, the phrase $h(x^k)^{-1}e'{}^p_{2i+1}h(x^k) =
(e_1f_2)^{-k}e'{}^p_{2i+1}(e_1f_2)^{k}$  has a single amalgamation that
reduces the length by exactly two.  

\smallbreak
Now upper bounds for the length of a word $w$ are computed which depend
on the variable $\sigma$:  

\medbreak\noindent
{\bf 3.15.1.a.}  The integer $\sigma$ is even.   Here
$\length(w) \geq 2kj + (j-1) + \epsilon_1 + \epsilon_{2j+1} - \sigma
-5(j-1- \sigma) \geq  2kj + (j-1) - \sigma -5(j-1- \sigma) -1 = 2kj -
4j +4\sigma + 3 \geq 2j(k-2) + 3$.

\medbreak\noindent
{\bf 3.15.1.b.}  The integer $\sigma$ is odd and there are more
occurences of the phrase $$h(x^k)e'{}^p_{2i+1}h(x^k)^{-1}$$
than of the phrase $h(x^k)^{-1}e'{}^p_{2i+1}h(x^k)$. Here
$\length(w) \geq 2kj + (j-1) + \epsilon_1 + \epsilon_{2j+1} - (\sigma -
1) -5(j-1- \sigma) 
\geq 2kj + (j-1) - 2 - (\sigma - 1) -5(j-1- \sigma)  = 2kj - 4j +4\sigma
+ 3 \geq 2j(k-2) + 3$.  

\medbreak\noindent
{\bf 3.15.1.c.}  The integer $\sigma$ is odd and there are more
occurences of the phrase

$h(x^k)^{-1}e'{}^p_{2i+1}h(x^k)$

\noindent  than of the phrase $h(x^k)e'{}^p_{2i+1}h(x^k)^{-1}$. Here
$\length(w) \geq 2kj + (j-1) + \epsilon_1 + \epsilon_{2j+1} - (\sigma -1)
- 2 -5(j-1- \sigma)
\geq 2kj + (j-1) - (\sigma + 1) -5(j-1- \sigma)  = 2kj - 4j +4\sigma + 3
\geq 2j(k-2) + 3$.

Note that in all cases the elements $e'$ and $h(x^k)$ cannot generate
any elements of  length less than $2k-1$ other than those elements of
form $e'{}^p$.  For instance, in 3.15.1.a or in  3.15.1.b, then the length
of the word $w$ is at least $2j(k-2) + 3$. Note if j = 0, the  resulting
word is of form $e'{}^p$. So $j \geq 1$ and if $k$ is at least three, then
the  inequality $2j(k-2) + 3< 2k-1$ implies that $j < 1$, a
contradicition.  But if $k=2$, then the inequality reduces to $3<3$,
another contraction.  The statement of 3.14 follows. 

Now if $h(bc)=f'$, repeating the above argument, replacing $x^k$ by
$x^{-k}$ and $e'$  by $f'$ yields the statement of 3.14.

\medbreak\noindent
{\bf 3.15.2.} The element $h(bc)$ lies in a conjugate of a free factor,
i.e the element $h(bc)$ is a conjugate of an element in a free factor of ${\bf Z}_p *{\bf Z}_q$.  

Without loss of generality suppose $h(bc) = ge'g^{-1}$, where the
element $g$ has one of the four normal forms listed below:

i:  $g = e_1' \cdots e_n'$.   

ii:  $g = e_1' \cdots f_n'$.   

iii:  $g = f_1' \cdots e_n'$.   

iv:  $g = f_1' \cdots f_n'$.

\noindent Note that in cases i and iii, the length of the element $h(bc)$ is
$2n - 1$ and in cases ii and iv, that the length of $h(bc)$ is $2n + 1$.  

Let $w$ be a word in the elements $h(bc)$ and $h(x^k)$.  Recall that
$h(x^k)$ has normal form $(e_1f_2)^k$.  The normal form for the element
$h(x^k)$ can be written out and  depends upon which of the four normal
forms above the element $g$ has.  Since ${\bf Z}_p$ is an abelian group,
we may assume for the purposes of working with $h(bc)$ that $g$ has
normal form of type ii or type iv.  If not the trivial element, the
element $h(bc)^p$ has length at least three.  By conjugating, if  necessary,
we may assume that the element $e_1'$ beginning the word
$g$ in case ii is not the inverse of the element $e_1$ that begins
the word $(e_1f_2)^k$.  Similarly, we may  assume that the element $f_1'$
beginning the word $g$ in case iv is not the inverse of the 
element $f_2$, ending the word $(e_1f_2)^k$.   The word $w$ has the
following form:

$w = h(bc)^{p_1}(e_1f_2)^{p_2 \cdot k}\cdots (e_1f_2)^{p_{2j} \cdot 
k}h(bc)^{p_{2j+1}}$.

\noindent With the possible exception of $p_1$ and $p_{2j+1}$, all
integers $p_i$ are nonzero and no power of $h(bc)$ is trivial.  For each
of cases ii and iv,  the phrase $h(bc)^p(e_1f_2)^{q \cdot  k}$ is
examined for cancellation and or amalgamation:

\smallbreak\noindent
{\bf 3.15.2.1.}  If $p$ is a positive integer,  writing out
the phrase using $(ge'g^{-1})^p$ instead of $h(bc)^p$ in case ii shows that
there is an amalgamation in the elements $e_1'{}^{-1}e_1$, and no
cancellation since the elements $e_1'$ and $e_1$ are not inverses.   In
case iv, there is no cancellation or amalgamation of any kind.  

If $p$ is a negative integer, there is no amalgamation or cancellation
in case ii.  In case iv, there is amalgamation, but no cancellation since the
elements $f_1'$ and $f_2$ are not inverse.  

\smallbreak\noindent
{\bf 3.15.2.2.}  For each of the cases, writing the phrase
$(e_1f_2)^kh(bc)^p(e_1f_2)^k$ shows that there is a single amalgamation
that reduces the uncancelled length by one.

\smallbreak\noindent
{\bf 3.15.2.3.}  For each of the cases, writing the phrase
$((e_1f_2)^{k})^{-1}h(bc)^p((e_1f_2)^{k})^{-1}$ shows that there is a
single  amalgamation that reduces the uncancelled length by one.

\smallbreak\noindent
{\bf 3.15.2.4.}  Writing the phrase
$(e_1f_2)^kh(bc)^p((e_1f_2)^{k})^{-1}$ shows that in case ii there is
no cancellation or amalgamation.  In case iv, there are two
amalgamations that reduce the uncancelled length by two.

\smallbreak\noindent
{\bf 3.15.2.5.}  Writing the phrase
$((e_1f_2)^{k})^{-1}h(bc)^p(e_1f_2)^k$ shows that in case ii,
there are two amalgamations that reduce the uncancelled length by two. 
In case vi there is no cancellation or amalgamation.

Let $w$ be an arbitary word of form $w = h(bc)^{p_1}(e_1f_2)^{p_2 \cdot
k}\cdots (e_1f_2)^{p_{2j} \cdot k}h(bc)^{p_{2j+1}}$.  For each of the
possible forms for $g$, an upper bound for the length of $w$ will now
be computed.

\smallbreak\noindent
{\bf 3.15.2.6.}  Assume that $g$ has type ii.  The length of
$w$ is computed  depending on whether $\sigma$ is even or odd.  

If $\sigma$ is even, then $\length(w)\geq 2kj + (j-1)(2n+1) - \sigma
-(j-1-\sigma) \geq (2n)(j-1)+2kj$.

If $\sigma$ is odd and there are more occurences of the phrase
$(e_1f_2)^kh(bc)^p((e_1f_2)^k)^{-1}$ in $w$ than of the phrase
$((e_1f_2)^k)^{-1}h(bc)^p(e_1f_2)^k$, then:

$\length(w)\geq 2kj + (j-1)(2n+1) -
(\sigma-1)-(j-1-\sigma)=2n(j-1) + 2kj +1$.
   
If $\sigma$ is odd and there are more occurences of the phrase
$((e_1f_2)^{k})^{-1}h(bc)^p(e_1f_2)^k$ in $w$ than of the phrase
$(e_1f_2)^kh(bc)^p((e_1f_2)^{k})^{-1}$, then:

$\length(w)\geq 2kj + (j-1)(2n+1) - (\sigma-1)-2
-(j-1-\sigma)=2kj + 2n(j-1) - 1$.

\smallbreak\noindent
{\bf 3.15.2.7.}  Assume that $g$ has type iv.  The length of
$w$ is computed  depending on whether $\sigma$ is even or odd.  

If $\sigma$ is even, then $\length(w)\geq 2kj + (j-1)(2n+1)- \sigma
-(j-1-\sigma) \geq (2n)(j-1)+2kj \geq 2kj = 2kj + 2n(j-1)$.

If $\sigma$ is odd and there are more occurences of the phrase
$(e_1f_2)^kh(bc)^p((e_1f_2)^{k})^{-1}$ in $w$ than of the phrase
$((e_1f_2)^{k})^{-1}h(bc)^p(e_1f_2)^k$, then:

$\length(w)\geq 2kj + (j-1)(2n+1) - (\sigma-1)-2
-(j-1-\sigma)=2kj +2n(j-1)-1$.
   
If $\sigma$ is odd and there are more occurences of the phrase
$((e_1f_2)^{k})^{-1}h(bc)^p(e_1f_2)^k$ in $w$ than of the phrase
$(e_1f_2)^kh(bc)^p((e_1f_2)^{k})^{-1}$, then:

$\length(w)\geq 2kj + (j-1)(2n+1) - (\sigma-1)
-(j-1-\sigma)=2kj +2n(j-1)+1$.
   
Recall that the integer $k$ is at least two and that the integer $n$ is at
least one.  Solving each of the following equations:

$2kj + 2n(j-1) - 1 < 2k -1$

$2kj + 2n(j-1)  < 2k -1$

$2kj + 2n(j-1) + 1 < 2k -1$

\noindent for the variable $j$, shows that if $w$ is the representative of
a word of length less than $2k-1$, then $w$ has form $h(bc)^p$.  But by
3.15.2,  such an element has odd length.  As the length of the element
$(e_1f_2)^j$ is the even number $2j$,  the elements $h(bc)$
and $h(x^k) = (e_1f_2)^k$ cannot generate an element of form
$(e_1f_2)^j$ for $1 \leq j \leq k-1$.  Theorem 3.14, and hence Theorem
3.10 follows.

\medbreak\noindent
{\bf 3.16.} ({\it Proof of 3.1.})  Theorem 3.1 now follows from 3.9 and 3.10. 
\hfill $\diamondsuit$

\end{document}